\documentclass[11pt]{article}
\usepackage{amsfonts}
\usepackage{latexsym,amsmath}
\usepackage{amssymb,array}
\usepackage[numbers]{natbib}
\usepackage{setspace}
\usepackage{enumitem}
\usepackage{xcolor}
\usepackage{graphicx}
\setlist{noitemsep}
\makeatletter \oddsidemargin 0in \evensidemargin 0in \textwidth 16cm
\begin{document}
    \newcommand{\ov}{\overline}
    \newcommand{\om}{\omega}
    \newcommand{\ga}{\gamma}
    \newcommand{\cd}{\circledast}
    \newtheorem{thm}{Theorem}[section]
    \newtheorem{remark}{Remark}[section]
    \newtheorem{counterexample}{Counterexample}[section]
    \newtheorem{coro}{Corollary}[section]
    \newtheorem{propo}{Proposition}[section]
    \newtheorem{definition}{Definition}[section]
    \newtheorem{example}{Example}[section]
    \newtheorem{lem}{Lemma}[section]
    \date{}

\title{Stochastic comparisons of lifetimes of series and parallel systems with dependent and heterogeneous components}
\author{Arindam Panja$^1$, Pradip Kundu$^2$\footnote
    {Corresponding author, e-mail: kundu.maths@gmail.com}
     and Biswabrata Pradhan$^1$\\
$^1$ SQC \& OR Unit, Indian Statistical Institute, Kolkata-700108, India\\
$^2$ Decision Science and Operations Management, Birla Global
University,\\ Bhubaneswar, Odisha-751003,
India\let\thefootnote\relax\footnote{Email Address:
arindampnj@gmail.com (Arindam Panja), bis@isical.ac.in (Biswabrata
Pradhan)}} \maketitle 
\begin{abstract}
This work considers stochastic comparisons of lifetimes of series
and parallel systems with dependent and heterogeneous components
having lifetimes following the proportional odds (PO) model. The
joint distribution of component lifetimes is modeled by Archimedean
survival copula. We discuss some potential applications of our
findings in system reliability and actuarial science.
\end{abstract}
\section{Introduction}\label{intro}Suppose $X_{1},X_{2},\ldots,X_{n}$ are the random variables denoting
the lifetimes of the components of a system with $n$ components.
Then the system lifetime is a function of $X_{1},X_{2},\ldots,X_{n}$.
Let $X_{k:n}$, $k=1,2,\ldots,n$ denote the $k$th order statistic
corresponding to the random variables $X_{1},X_{2},\ldots,X_{n}$.
Then the smallest and the largest order statistics $X_{1:n}$ and
$X_{n:n}$, respectively, represent the lifetimes of the series and
the parallel systems. There have been a number of works on
stochastic comparisons of system lifetimes where component lifetimes
follow different family of distributions
\cite{barm,ding,fangl,fangL,fangR,gupta,hazra,lif,nava}. However,
most of the works have considered mutual independence among the
concerned random variables. Recently, Fang et al. \cite{fangR}, Li
and Fang \cite{lif} and Li and Li \cite{lil} have considered
stochastic comparison of system lifetimes with dependent and
heterogeneous component lifetimes following the proportional hazard
rate (PHR) model.

The proportional odds (PO) model, as introduced by Bennet
\cite{benn} and later discussed by Kirmani and Gupta \cite{kirmani}
is a very important model in reliability theory and survival
analysis. Let $X$ and $Y$ be two random variables with distribution
functions $F_X(\cdot)$, $F_Y(\cdot)$, survival functions
$\bar{F}_X(\cdot)$, $\bar{F}_Y(\cdot)$ and hazard rate functions
$r_{X}(\cdot)$, $r_{Y}(\cdot)$ respectively. Let the odds functions
of $X$ and $Y$ be defined by $\tau_{X}(t)=\bar{F}_X(t)/F_X(t)$ and
$\tau_{Y}(t)=\bar{F}_Y(t)/F_Y(t)$, respectively. If the random
variable $X$ represents the lifetime of a component, then the odds
function $\tau_{X}(t)$ represents the odds of that component
functioning beyond time $t$. The random variables $X$ and $Y$ are
said to satisfy PO model if $\tau_{Y}(t)=\alpha \tau_{X}(t)$ for all
admissible $t$, where $\alpha$ is a proportionality constant known
as proportional odds ratio. Then the survival functions of $X$ and
$Y$ are related as
\begin{equation}\label{poalt}\bar{F}_Y(t)=\frac{\alpha
\bar{F}_X(t)}{1-\bar{\alpha}\bar{F}_X(t)},\end{equation} where
$\bar{\alpha}=1-\alpha$. We will say that the random variable $Y$ is
following the PO model with baseline survival function
$\bar{F}_X(\cdot)$ and parameter (proportionality constant)
$\alpha$. For easy interpretation, we can think of $X$ as the
lifetime of a member of control group, and $Y$ as that of a member
of treatment group. For two random variables satisfying the PO
model, the ratio of hazard rates converges to unity as time tends to
infinity, which is in contrast to the PHR model where this ratio
remains constant with time. The convergence property of hazard
functions makes the PO model reasonable in many practical
applications as discussed in \cite{benn,kirmani,ross}. For more
applications of PO model one may refer to \cite{coll,zha}. Also, the
model (\ref{poalt}), with $0<\alpha <\infty$, provides us a method
of generating more flexible new family of distributions by
introducing the parameter $\alpha$ to an existing family of
distributions \cite{marsh1}. The family of distributions so obtained
is known as Marshall-Olkin family of distributions
\cite{cord1,marsh1}. Thus, model (\ref{poalt}) has implications both
in terms of the PO model and in generating new family of flexible
distributions, which makes it worth investigating.

Let $\textbf{X}=(X_{1},X_{2},\ldots,X_{n})$ be a random vector with
joint distribution function $F(\cdot)$ and joint survival function
$\bar{F}(\cdot)$. Also let the distribution function and the
survival function of $X_i$ are $F_i(\cdot)$ and $\bar{F}_i(\cdot)$
respectively for $i=1,2,\ldots,n$. The joint distribution of
$X_{1},X_{2},\ldots,X_{n}$ can be represented by a copula model. If
there exist $K:[0,1]^n\mapsto [0,1]$ and $\bar{K}:[0,1]^n\mapsto
[0,1]$ such that $F(x_1,\ldots,x_n)=K(F_1(x_1),\ldots,F_n(x_n))$ and
$\bar{F}(x_1,\ldots,x_n)=\bar{K}(\bar{F}_1(x_1),\ldots,\bar{F}_n(x_n))$
for all $x_i$, $i\in I_n$, then $K$ and $\bar{K}$ are called the
copula and survival copula of $X$, respectively. If
$\varphi:[0,+\infty)\mapsto [0,1]$ with $\varphi(0)=1$ and
$\lim_{t\rightarrow +\infty} \varphi(t)=0$ is $(n-2)$th
differentiable, then
$K_{\varphi}(u_1,\ldots,u_n)=\varphi(\varphi^{-1}(u_1)+\ldots+\varphi^{-1}(u_n))=\varphi(\sum_{i=1}^n
\phi (u_i))$ for all $u_i\in(0,1]$, $i\in I_n$ is called an
Archimedean survival copula with generator $\varphi$ provided
$(-1)^{k}\varphi^{(k)}(t)\geq 0$, $k=0,1,\ldots,n-2$ and
$(-1)^{n-2}\varphi^{(n-2)}(t)$ is decreasing and convex for all
$t\geq 0$. Here $\phi=\varphi^{-1}$ is the right continuous inverse
of $\varphi$ so that $\phi(u)=\varphi^{-1}(u)=\sup\{t\in
\mathbb{R}:\varphi(t)>u\}$. Navarro and Spizzichino \cite{nava} have
derived usual stochastic ordering for lifetimes of series and
parallel systems having component lifetimes sharing a common copula,
with the idea of mean reliability function associated with the
common copula. Li and Fang \cite{lif} investigated stochastic order
between two samples of dependent random variables following PHR
model and having Archimedean survival copula. Fang et al.
\cite{fangR} derived some stochastic ordering results for minimum as
well as for maximum of samples equipped with Archimedean survival
copulas and following PHR model and proportional reversed
hazard rate (PRH) model, respectively. 
Li and Li \cite{lil} investigated hazard rate order on minimums of
sample following PHR model, and reversed hazard rate order on
maximums of sample following PRH model, where both the samples
coupled with Archimedean survival copula.

In case of PO model, some authors, e.g. Kundu and Nanda
\cite{kundu1}, Kundu et al. \cite{kundu2}, Nanda and Das
\cite{nanda} have investigated stochastic comparison of systems with
independent components. However, to the best of our knowledge, no
research work has been done on stochastic comparison of system
lifetimes with dependent and heterogeneous component lifetimes
following PO model. In this work, we investigate stochastic
comparisons of lifetimes of series and parallel systems with
dependent and heterogenous components having lifetimes following the
PO model. The joint distribution of component lifetimes is modeled
by Archemedian survival copula. It is shown that the usual
stochastic ordering and hazard rate ordering hold for series systems
under certain conditions whereas for parallel system stochastic
ordering and reversed hazard rate ordering hold.
\par The organization of the paper is as follows. Section 2
recalls some definitions of majorization, stochastic orders, and
some lemmas used in the sequel. In Section 3, we investigate
stochastic comparisons between series systems of dependent and
heterogenous components having lifetimes following the PO model and
coupled by Archimedean survival copulas. Section 4 investigates the
same in case of parallel systems. Section 5 presents some potential
applications of the proposed results.
\section{Preliminaries}\label{preli}
Given a vector $\mbox{\boldmath $x$}=(x_1,x_2,...,x_n)\in
\mathbb{R}^n$, denote $x_{(1)}\leq x_{(2)}\leq ...\leq x_{(n)}$ as
increasing arrangement of $x_1,x_2,\ldots,x_n$.
\begin{definition}
Let $\mbox{\boldmath $x$}=(x_1,x_2,\dots,x_n)$ and $\mbox{\boldmath
$y$}=(y_1,y_2,\dots,y_n)$ in $\mathbb{R}^n$ be any two vectors.
\begin{enumerate}
\item[(i)] The vector $\mbox{\boldmath $x$}$ is said to majorize the vector $\mbox{\boldmath $y$}$, i.e., $\mbox{\boldmath $x$}$ is larger than $\mbox{\boldmath $y$}$ in majorization order (denoted as $\mbox{\boldmath $x$}\stackrel{m}{\succeq}\mbox{\boldmath $y$}$)
if (cf. \cite{marsh3})
\begin{equation*}
\sum_{i=1}^j x_{(i)}\le\sum_{i=1}^j y_{(i)},~\text{for
all}\;j=1,\;2,\;\ldots, n-1,\;\;and \;\;\sum_{i=1}^n
x_{(i)}=\sum_{i=1}^n y_{(i)}.
\end{equation*}
\item [(ii)] The vector $\mbox{\boldmath $x$}$ is said to weakly supermajorize the vector $\mbox{\boldmath
$y$}$, denoted as $\mbox{\boldmath $x$}\stackrel{ w}{\succeq}
\mbox{\boldmath $y$}$ if (cf. \cite{marsh3})
 \begin{eqnarray*}
  \sum\limits_{i=1}^j x_{(i)}\leq \sum\limits_{i=1}^j y_{(i)},~\text{for all}\;j=1,2,\dots,n.
 \end{eqnarray*}
 \item [(iii)] The vector $\mbox{\boldmath $x$}$ is said to be $p$-larger than the vector $\mbox{\boldmath $y$}$
 (denoted as $\mbox{\boldmath $x$}\stackrel{ p}{\succeq} \mbox{\boldmath $y$}$)
 if (cf. \cite{bon1})
 \begin{eqnarray*}
  \prod\limits_{i=1}^j x_{(i)}\leq \prod\limits_{i=1}^j y_{(i)},~ \text{for all}\;j=1,2,\dots,n.
 \end{eqnarray*}
\end{enumerate}
\end{definition}
It can be seen that
$$\mbox{\boldmath x}\stackrel{m}{\succeq}\mbox{\boldmath
y}\Rightarrow\mbox{\boldmath x}\stackrel{ w}{\succeq}
\mathbf{y}\Rightarrow\mbox{\boldmath x}\stackrel{ p}{\succeq}
\mathbf{y}.\hfill\Box $$
\begin{definition} \cite{shaked} Let $X$ and $Y$ be nonnegative absolutely
continuous random variables with cumulative distribution
functions $F_X(\cdot)$, $F_Y(\cdot)$, survival functions
$\bar{F}_X(\cdot)$, $\bar{F}_Y(\cdot)$, hazard (failure) rate
functions $r_X(\cdot)$, $r_Y(\cdot)$, and the reversed hazard rate
functions $\tilde r_X(\cdot)$ and $\tilde r_Y(\cdot)$, respectively.
Then $X$ is said to be smaller than $Y$ in the
\begin{enumerate}[label=(\roman*)]
\item usual stochastic order (denoted as $X
\leq_{st} Y$) if $\bar{F}_X(t)\leq \bar{F}_Y(t)$ for all $t$;
\item hazard rate order (denoted as $X \leq_{hr} Y$) if
$\bar{F}_Y(t)/\bar{F}_X(t)$ is increasing in $t\geq 0$, or
equivalently if $r_X(t)\geq r_Y(t)$ for all $t\geq 0$;
\item reversed hazard rate order (denoted as $X \leq_{rhr} Y$) if
$F_Y(t)/F_X(t)$ is increasing in $t>0$, or equivalently if $\tilde
r_X(t)\leq \tilde r_Y(t)$ for all $t>0$. $\hfill\Box$
\end{enumerate}
\end{definition}
\begin{lem}\label{lesc} \cite{marsh3}
Let $I \subseteq \mathbb{R}$ be an open interval and let $\zeta:
I^n\rightarrow \mathbb{R}$ be continuously differentiable. Necessary
and sufficient conditions for $\zeta$ to be Schur-convex (resp.
Schur-concave) on $I^n$ are that $\zeta$ is symmetric on $I^n$, and
for all $i\neq j$,
$$(u_i-u_j)\left(\zeta_{(i)}({\mbox{\boldmath
$u$}})-\\ \zeta_{(j)}({\mbox{\boldmath $u$}})\right)\geq
(\text{resp.}\leq)\;0~\text{for all}~{\mbox{\boldmath
$u$}}=(u_1,u_2,...,u_n)\in I^n,$$ where
$\zeta_{(k)}({\mbox{\boldmath $u$}})=\partial\zeta({\mbox{\boldmath
$u$}})/\partial u_k$.
 \end {lem}
 \begin{lem}\label{lewm} \cite{marsh3}
 Let $\mathcal{A} \subseteq \mathbb{R}^n$, and $\zeta: \mathcal{A}\rightarrow \mathbb{R}$ be a function. Then, for ${\mbox{\boldmath $x$}},{\mbox{\boldmath $y$}}\in \mathcal{A},$
$${\mbox{\boldmath $x$}}\stackrel{w}\succeq {\mbox{\boldmath $y$}}\; \implies\;\zeta({\mbox{\boldmath $x$}})\geq (\text{resp. }\leq)\; \zeta({\mbox{\boldmath $y$}})$$
if and only if $\zeta$ is both decreasing (resp. increasing) and
Schur-convex (resp. Schur-concave) on $\mathcal{A}$.
 \end {lem}
\begin{lem}\label{lesce} \cite{khale}
  Let $\zeta :(0,\infty)^n\rightarrow \mathbb{R}$ be a function.
 Then,
 $${\mbox{\boldmath $x$}}\stackrel{\rm p}\succeq {\mbox{\boldmath $y$}}\implies \zeta({\mbox{\boldmath $x$}})\geq(\text{resp.}\leq) \;\zeta({\mbox{\boldmath $y$}})$$
if and only if the following two conditions hold:
\begin{enumerate}
 \item[(i)]$\zeta(e^{v_1},\dots,e^{v_n})$ is Schur-convex (resp. Schur-concave) in $(v_1,\dots,v_n)$,
\item[(ii)]$\zeta(e^{v_1},\dots,e^{v_n})$ is decreasing (resp. increasing) in each $v_i,$ for $i=1,\dots,n,$
\end{enumerate}
where $v_i=\ln x_i$, for $i=1,\dots,n.$ $\hfill\Box$
\end{lem}
\begin{lem}\label{lecop} \cite{fangR}
For two $n$-dimensional Archimedean copulas $K_{\varphi_1}$ and
$K_{\varphi_2}$, if $\phi_2\circ\varphi_1$ is superadditive, then
$K_{\varphi_1}(\textbf{u})\leq K_{\varphi_2}(\textbf{u})$ for all
$\textbf{u}\in [0,1]^n$.
\end {lem}
\section{Series systems with dependent and heterogeneous component lifetimes following PO Model}
Here, we consider the comparison of lifetimes of two series systems
with heterogeneous and dependent components. We assume that the
lifetime vector $X=(X_1,X_2,...,X_n)$ is a set of dependent random
variables coupled with Archimedean survival copula with generator
$\varphi$ and following the PO model with baseline survival function
$\bar{F}$, denoted as $X\sim PO(\bar{F},\boldsymbol\alpha,\varphi)$,
where $\boldsymbol\alpha=(\alpha_1,\alpha_2,...,\alpha_n)\in
\mathbb{R}^n_{+}$ is the proportional odds ratio vector. The
survival function and the hazard rate function of $X_i$ are
$$\bar{F}_{\alpha_i}(x)=\frac{\alpha_i
\bar{F}(x)}{1-\bar{\alpha}_i\bar{F}(x)}
~\text{and}~r_{\alpha_i}(x)=\frac{r(x)}{1-\bar{\alpha}_i\bar{F}(x)},
~\text{respectively},$$ where $\bar{\alpha}_i=1-\alpha_i$,
$i=1,\ldots,n$ and $r$ denotes the baseline hazard rate function.
The survival functions of $X_{1:n}$ is given by
\begin{equation}\label{ssf}\bar{F}_{X_{1:n}}(x)=P(X_{1:n}>x)=P(X_{i}>x,i\in
I_n)=\varphi\left(\sum_{i=1}^n
\phi\left(\bar{F}_{\alpha_i}(x)\right)\right)=S_1(\bar{F}(x),\boldsymbol\alpha,\varphi),~say,\end{equation}
where $\phi(u)=\varphi^{-1}(u)$, $u\in(0,1]$.\\ The hazard rate
function of $X_{1:n}$ is obtained as
\begin{equation}\label{shr}r_{X_{1:n}}(x)=r(x)\frac{\varphi'\left(\sum_{i=1}^n
\phi\left(\bar{F}_{\alpha_i}(x)\right)\right)}{\varphi\left(\sum_{i=1}^n
\phi\left(\bar{F}_{\alpha_i}(x)\right)\right)}\sum_{i=1}^n
\phi'\left(\bar{F}_{\alpha_i}(x)\right)\frac{\bar{F}_{\alpha_i}(x)}{1-\bar{\alpha}_i\bar{F}(x)}.\end{equation}
\begin{lem}\label{lej1}
For any $x\in[0,1]$, $S_1(x,\boldsymbol\alpha,\varphi)$ is
increasing in $\alpha_i$, $i\in I_n$. Furthermore $S_1$ is
Schur-concave with respect to $\boldsymbol\alpha$.
 \end{lem}
\textbf{Proof:} For $s\in I_n$,
$$\frac{\partial S_1}{\partial \alpha_s}=\varphi'\left(\sum_{i=1}^n
\phi\left(\frac{\alpha_i x}{1-\bar{\alpha}_i
x}\right)\right)\phi'\left(\frac{\alpha_s x}{1-\bar{\alpha}_s
x}\right)\frac{x(1-x)}{(1-\bar{\alpha}_s x)^2}.$$ Since both
$\varphi(u)$ and $\phi(u)$ are decreasing for all $u\geq 0$,
$\frac{\partial S_1}{\partial \alpha_s}\geq0$. As a result
$S_1(x,\boldsymbol\alpha,\varphi)$ is increasing in $\alpha_i$,
$i\in I_n$ for any $x\in[0,1]$.\\For $s\neq t$,
\begin{eqnarray}\nonumber &&(\alpha_s-\alpha_t)\left(\frac{\partial S_1}{\partial
\alpha_s}-\frac{\partial S_1}{\partial \alpha_t}\right)\\\nonumber
&=&(\alpha_s-\alpha_t)\varphi'\left(\sum_{i=1}^n
\phi\left(\frac{\alpha_i x}{1-\bar{\alpha}_i
x}\right)\right)\left[\phi'\left(\frac{\alpha_s x}{1-\bar{\alpha}_s
x}\right)\frac{x(1-x)}{(1-\bar{\alpha}_s
x)^2}-\phi'\left(\frac{\alpha_t x}{1-\bar{\alpha}_t
x}\right)\frac{x(1-x)}{(1-\bar{\alpha}_t x)^2}\right]\\
\label{eqlem}
&\stackrel{sign}{=}&(\alpha_s-\alpha_t)\left(-\varphi'\left(\sum_{i=1}^n
\phi\left(u_i\right)\right)\right)\left[\left(-\phi'\left(u_s\right)\right)\frac{1}{(1-\bar{\alpha}_s
x)^2}-\left(-\phi'\left(u_t\right)\right)\frac{1}{(1-\bar{\alpha}_t
x)^2}\right],
\end{eqnarray}where $u_i=\frac{\alpha_i x}{1-\bar{\alpha}_i x}$ and $\stackrel{`sign'}{=}$ means equal in sign. Since both $\varphi$ and $\phi$ are decreasing, and $\phi'$ is increasing, it follows from (\ref{eqlem}) that $(\alpha_s-\alpha_t)\left(\frac{\partial S_1}{\partial
\alpha_s}-\frac{\partial S_1}{\partial \alpha_t}\right)\leq 0$. So
from Lemma \ref{lesc}, $S_1$ is Schur-concave in
$\boldsymbol\alpha=(\alpha_1,\alpha_2,...,\alpha_n)$.$\hfill\Box$

Suppose there are two series systems formed out of $n$ statistically
dependent and heterogeneous components where the component lifetimes
follow the PO model. The joint distribution of lifetimes of
components is represented by Archimedean copula. Consider two such
series systems with lifetime vectors $X=(X_1,X_2,...,X_n)$ and
$Y=(Y_1,Y_2,...,Y_n)$ having respective proportionality odds ratio
vectors $\boldsymbol\alpha=(\alpha_1,\alpha_2,...,\alpha_n)$ and
$\boldsymbol\beta=(\beta_1,\beta_2,...,\beta_n)$, where
$\boldsymbol\alpha, \boldsymbol\beta \in \mathbb{R}^n_{+}$. The
following theorem compares the lifetimes of these series systems in
the sense of usual stochastic order.
\begin{thm}\label{thst}Suppose the lifetime vectors $X\sim PO(\bar{F},\boldsymbol\alpha,\varphi_1)$ and $Y\sim
PO(\bar{F},\boldsymbol\beta,\varphi_2)$. If $\varphi_1$ or
$\varphi_2$ is log-convex and $\phi_2\circ\varphi_1$ is
superadditive, then
$$\boldsymbol\alpha \stackrel{p}{\succeq}\boldsymbol\beta~\text{implies}~ X_{1:n}\leq_{st}Y_{1:n}.$$
\end{thm}
\textbf{Proof:} Write $v_i=\ln \alpha_i$, $i=1,2,...,n$. Then as per
(\ref{ssf}),
$$\bar{F}_{X_{1:n}}(x)=\varphi_1\left(\sum_{i=1}^n
\phi_1\left(\frac{e^{v_i}
\bar{F}(x)}{1-(1-e^{v_i})\bar{F}(x)}\right)\right)=S_1(\bar{F}(x),(e^{v_1},e^{v_2},...,e^{v_n}),\varphi_1).$$
Here $S_1(\bar{F}(x),(e^{v_1},e^{v_2},...,e^{v_n}),\varphi_1)$ is
symmetric with respect to $(v_1,v_2,...,v_n)\in \mathbb{R}^n$. Now,
for $s\in I_n$,
$$\frac{\partial S_1}{\partial v_s}=\varphi_1^{'}\left(\sum_{i=1}^n
\phi_1\left(\frac{e^{v_i} x}{1-(1-e^{v_i})
x}\right)\right)\phi_1^{'}\left(\frac{e^{v_s} x}{1-(1-e^{v_s})
x}\right)\frac{x(1-x)e^{v_s}}{(1-(1-e^{v_s}) x)^2},$$ so that
$S_1(x,(e^{v_1},e^{v_2},...,e^{v_n}),\varphi_1)$ is increasing in
each $v_i$, $i=1,2,...,n$ for any $x\in[0,1]$.\\ Now, for $s\neq t$,
\begin{eqnarray}\nonumber &&(v_s-v_t)\left(\frac{\partial S_1}{\partial
v_s}-\frac{\partial S_1}{\partial v_t}\right)\\
\nonumber &=&(v_s-v_t)\left(-\varphi_1^{'}\left(\sum_{i=1}^n
\phi_1\left(\frac{e^{v_i} x}{1-(1-e^{v_i})
x}\right)\right)\right)\left[\left(-\phi_1^{'}\left(\frac{e^{v_s}
x}{1-(1-e^{v_s}) x}\right)\right)\frac{xe^{v_s}}{(1-(1-e^{v_s})
x)^2}\right.\\\nonumber
&&\left.-\left(-\phi_1^{'}\left(\frac{e^{v_t} x}{1-(1-e^{v_t})
x}\right)\right)\frac{xe^{v_t}}{(1-(1-e^{v_t}) x)^2}\right]\\
\label{eqth3.1}
&\stackrel{sign}{=}&(v_s-v_t)\left[\left(-\frac{\varphi_1(\phi_1(u_s))}{\varphi_1^{'}(\phi_1(u_s))}\right)\frac{1
}{1-(1-e^{v_s})
x}-\left(-\frac{\varphi_1(\phi_1(u_t))}{\varphi_1^{'}(\phi_1(u_t))}\right)\frac{1
}{1-(1-e^{v_t}) x}\right],
\end{eqnarray}
where $u_s=\frac{e^{v_s} x}{1-(1-e^{v_s})x}$. If $\varphi_1$ is
log-convex, from (\ref{eqth3.1}) it follows that
$(v_s-v_t)\left(\frac{\partial S_1}{\partial v_s}-\frac{\partial
S_1}{\partial v_t}\right)\leq 0$. Hence from Lemma \ref{lesc},
$S_1(x,(e^{v_1},e^{v_2},...,e^{v_n}),\varphi_1)$ is Schur-concave in
$(v_1,v_2,...,v_n)$ if $\varphi_1$ is log-convex. Then from Lemma
\ref{lesce}, we have
\begin{equation}\label{eqth3.11}\boldsymbol\alpha
\stackrel{p}{\succeq}\boldsymbol\beta~\text{implies}~S_1(\bar{F}(x),\boldsymbol\alpha,\varphi_1)\leq
S_1(\bar{F}(x),\boldsymbol\beta,\varphi_1).\end{equation} Since
$\phi_2\circ\varphi_1$ is superadditive, from Lemma \ref{lecop}, we
have
\begin{equation}\label{eqth3.12}S_1(\bar{F}(x),\boldsymbol\beta,\varphi_1)\leq
S_1(\bar{F}(x),\boldsymbol\beta,\varphi_2).\end{equation} Thus
combining (\ref{eqth3.11}) and (\ref{eqth3.12}) we get
$S_1(\bar{F}(x),\boldsymbol\alpha,\varphi_1)\leq
S_1(\bar{F}(x),\boldsymbol\beta,\varphi_2)$, that is
$X_{1:n}\leq_{st}Y_{1:n}$.\\ Now suppose $\varphi_2$ is log-convex,
then $S_1(\bar{F}(x),\boldsymbol\alpha,\varphi_2)\leq
S_1(\bar{F}(x),\boldsymbol\beta,\varphi_2)$. Since
$\phi_2\circ\varphi_1$ is superadditive, we have
$S_1(\bar{F}(x),\boldsymbol\alpha,\varphi_1)\leq
S_1(\bar{F}(x),\boldsymbol\alpha,\varphi_2)$. So combining we get
$X_{1:n}\leq_{st}Y_{1:n}$.$\hfill\Box$
\begin{coro}\label{corst}Suppose the lifetime vectors $X\sim PO(\bar{F},\boldsymbol\alpha,\varphi)$ and $Y\sim
PO(\bar{F},\boldsymbol\beta,\varphi)$. If $\varphi$ is log-convex,
then
$$\boldsymbol\alpha \stackrel{p}{\succeq}\boldsymbol\beta~\text{implies}~
X_{1:n}\leq_{st}Y_{1:n}. ~~~\hfill\Box $$
\end{coro}
The following counterexample shows that one may not get the the
ordering result in Theorem \ref{thst} if the sufficient conditions
on the generator functions are dropped.
    \begin{counterexample}\label{exthst}\normalfont
        Consider two series systems, each
        comprising of three dependent and heterogeneous components with respective survival
        functions $\bar{F}_{X_{1:3}}(x)=S_1(\bar{F}(x),\boldsymbol\alpha,\varphi_1)$ and
        $\bar{F}_{Y_{1:3}}(x)=S_1(\bar{F}(x),\boldsymbol\beta,\varphi_2)$ with
        $\bar{F}(x)=e^{-(x)^{1.5}}$, $x\geq
        0$, $\boldsymbol\alpha=(2,3,5.5)$, $\boldsymbol\beta=(2.5,3.5,3.8)$ so that $\boldsymbol\alpha
        \stackrel{p}{\succeq}\boldsymbol\beta$. First we take  $\varphi_1(x)=(2/(1+e^x))^{1/\theta}$,
        $\theta=0.9$ and  $\varphi_2(x)=e^{1-(1+x)^{1/\eta}}$ with
        $\eta=0.3$ so that $\phi_2\circ \varphi_1$ is not super
        additive, and $\varphi_1$ and $\varphi_2$ are not
        log-convex. We depict $\bar{F}_{X_{1:3}}(x)$ and $\bar{F}_{Y_{1:3}}(x)$ in Figure \ref{picth3.1p} for some finite range of $x$.
        From this figure we observe that the stochastic ordering result in Theorem~\ref{thst} is not attained.$\hfill\Box$
        \begin{figure}\begin{center}
                \includegraphics[width=10cm]{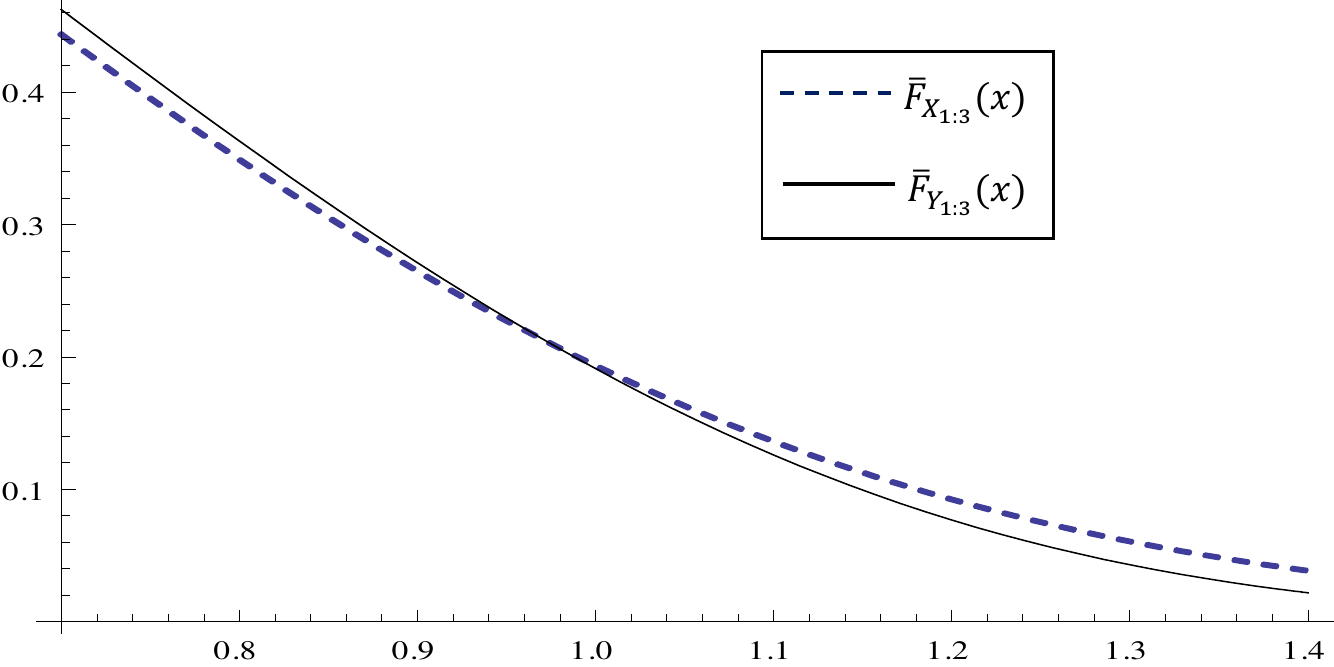}\\
                \caption{Plots of $\bar{F}_{X_{1:3}}(x)$ and $\bar{F}_{Y_{1:3}}(x)$ when $\phi_2\circ \varphi_1$ is not super additive, and $\varphi_1$ and $\varphi_2$ are not log-convex.}\label{picth3.1p}\end{center}
        \end{figure}
    \end{counterexample}
Since $p$-larger order is weaker than weakly supermajorization
order, the following theorem shows that we can get the ordering
result in Theorem \ref{thst} under weakly supermajorization order
with fewer condition.
\begin{thm}\label{thst1}Suppose the lifetime vectors $X\sim PO(\bar{F},\boldsymbol\alpha,\varphi_1)$ and $Y\sim
PO(\bar{F},\boldsymbol\beta,\varphi_2)$. If $\phi_2\circ\varphi_1$
is superadditive, then
$$\boldsymbol\alpha \stackrel{w}{\succeq}\boldsymbol\beta~\text{implies}~ X_{1:n}\leq_{st}Y_{1:n}.$$
\end{thm}
\textbf{Proof:} From Lemma \ref{lej1} and Lemma \ref{lewm}, we have
$$\boldsymbol\alpha
\stackrel{w}{\succeq}\boldsymbol\beta~\text{implies}~S_1(\bar{F}(x),\boldsymbol\alpha,\varphi_1)\leq
S_1(\bar{F}(x),\boldsymbol\beta,\varphi_1).$$ Since
$\phi_2\circ\varphi_1$ is superadditive, so from Lemma \ref{lecop},
we have, $S_1(\bar{F}(x),\boldsymbol\beta,\varphi_1)\leq
S_1(\bar{F}(x),\boldsymbol\beta,\varphi_2)$. Combining the above
results we have $S_1(\bar{F}(x),\boldsymbol\alpha,\varphi_1)\leq
S_1(\bar{F}(x),\boldsymbol\beta,\varphi_2)$. That is
$X_{1:n}\leq_{st}Y_{1:n}$.$\hfill\Box$
\begin{remark}\normalfont
It is to be noted that super-additive assumption of
$\phi_2\circ\varphi_1$ is satisfied by many members of Archimedean
survival copulas. For example, Archimedean survival copula with
generators (i) $\varphi_1(t)=e^{1-(1+t)^{\frac{1}{\theta}}}$ and
$\varphi_2(t)=\frac{\theta}{\log(e^{\theta}+t)}$, where $0<\theta
\leq 1$, (ii) $\varphi_1(t)=\frac{\theta}{\log(e^{\theta}+t)}$ and
$\varphi_2(t)=\log (e+t)^{-1/\theta}$, where $\theta>1$ and (iii)
$\varphi_1(t)=e^{1-(1+t)^{\frac{1}{\theta_1}}}$ and
$\varphi_2(t)=e^{1-(1+t)^{\frac{1}{\theta_2}}}$, where $\theta_2
\geq \theta_1 \geq 1$, satisfy super-additivity.
\end{remark}
\begin{coro}\label{corst1}Suppose the lifetime vectors $X\sim PO(\bar{F},\boldsymbol\alpha,\varphi)$ and $Y\sim
PO(\bar{F},\boldsymbol\beta,\varphi)$. Then
$$\boldsymbol\alpha \stackrel{w}{\succeq}\boldsymbol\beta~\text{implies}~
X_{1:n}\leq_{st}Y_{1:n}.$$
\end{coro}
\begin{lem}\label{lelu}
$I_1(\textbf{u})=\frac{\varphi'\left(\sum_{i=1}^n
u_i\right)}{\varphi\left(\sum_{i=1}^n u_i\right)}
\sum_{i=1}^n\frac{\varphi\left(u_i\right)}{\varphi'\left(u_i\right)}\left(1-\varphi\left(u_i\right)\right)$
is increasing in $u_s$, $s\in I_n$ and Schur-convex with respect to
$\textbf{u}=(u_1,...,u_n)$ if $\varphi$ is log-concave and
$\frac{\varphi (1-\varphi)}{\varphi'}$ is decreasing and concave.
 \end{lem}
\textbf{Proof:} Here $I_1(\textbf{u})$ is symmetric in $\textbf{u}$.
For $s\in I_n$,
$$\frac{\partial I_1(\textbf{u})}{\partial u_s}=\frac{\partial}{\partial u_s}\left(\frac{\varphi'\left(\sum_{i=1}^n
u_i\right)}{\varphi\left(\sum_{i=1}^n
u_i\right)}\right)\sum_{i=1}^n\frac{\varphi\left(u_i\right)}{\varphi'\left(u_i\right)}\left(1-\varphi\left(u_i\right)\right)+\frac{\varphi'\left(\sum_{i=1}^n
u_i\right)}{\varphi\left(\sum_{i=1}^n
u_i\right)}\frac{\partial}{\partial
u_s}\left(\frac{\varphi\left(u_s\right)}{\varphi'\left(u_s\right)}\left(1-\varphi\left(u_s\right)\right)\right).$$
Since $\varphi$ is log-concave, $\frac{\partial}{\partial
u_s}\left(\frac{\varphi'\left(\sum_{i=1}^n
u_i\right)}{\varphi\left(\sum_{i=1}^n u_i\right)}\right)\leq 0$. As
$\varphi (1-\varphi)/\varphi'$ is decreasing,
$\frac{\partial}{\partial
u_s}\left(\frac{\varphi\left(u_s\right)}{\varphi'\left(u_s\right)}\left(1-\varphi\left(u_s\right)\right)\right)\leq
0$. Then using the fact that $\varphi$ is deceasing, we have
$\frac{\partial I_1(\textbf{u})}{\partial u_s}\geq 0$. So
$I_1(\textbf{u})$ is increasing in $u_s$ for any $s\in I_n$. For
$s,t\in I_n$ with $s\neq t$,
$$\frac{\partial}{\partial
u_s}\left(\frac{\varphi'\left(\sum_{i=1}^n
u_i\right)}{\varphi\left(\sum_{i=1}^n
u_i\right)}\right)=\frac{\partial}{\partial
u_t}\left(\frac{\varphi'\left(\sum_{i=1}^n
u_i\right)}{\varphi\left(\sum_{i=1}^n u_i\right)}\right).$$ Then
\begin{eqnarray*}&&(u_s-u_t)\left(\frac{\partial I_1(\textbf{u})}{\partial
u_s}-\frac{\partial I_1(\textbf{u})}{\partial
u_t}\right)\\&=&(u_s-u_t)\frac{\varphi'\left(\sum_{i=1}^n
u_i\right)}{\varphi\left(\sum_{i=1}^n
u_i\right)}\left[\frac{\partial}{\partial
u_s}\left(\frac{\varphi\left(u_s\right)}{\varphi'\left(u_s\right)}\left(1-\varphi\left(u_s\right)\right)\right)-\frac{\partial}{\partial
u_t}\left(\frac{\varphi\left(u_t\right)}{\varphi'\left(u_t\right)}\left(1-\varphi\left(u_t\right)\right)\right)\right]\geq
0,\end{eqnarray*} where the inequality follows from the fact that
$\frac{\varphi (1-\varphi)}{\varphi'}$ is concave. So from lemma
\ref{lesc}, $I_1(\textbf{u})$ is Schur-convex with respect to
$\textbf{u}$. $\hfill\Box$

Next we show hazard rate ordering of two series systems formed out
of $n$ statistically dependent and heterogeneous components having
lifetimes following PO model.
\begin{thm}\label{thhr}Suppose the lifetime vectors $X\sim PO(\bar{F},\boldsymbol\alpha,\varphi)$ and $Y\sim
PO(\bar{F},\boldsymbol\beta,\varphi)$. If $\varphi$ is log-concave
and $\frac{\varphi (1-\varphi)}{\varphi'}$ is decreasing and concave
(or convex), then
$$\boldsymbol\alpha \stackrel{w}\succeq\boldsymbol\beta~\text{implies}~ X_{1:n}\leq_{hr}Y_{1:n}.$$
\end{thm}
\textbf{Proof:} From (\ref{shr}), we have
\begin{eqnarray*}r_{X_{1:n}}(x)&=&r(x)\frac{\varphi'\left(\sum_{i=1}^n
\phi\left(\bar{F}_{\alpha_i}(x)\right)\right)}{\varphi\left(\sum_{i=1}^n
\phi\left(\bar{F}_{\alpha_i}(x)\right)\right)}\sum_{i=1}^n
\phi'\left(\bar{F}_{\alpha_i}(x)\right)\frac{\bar{F}_{\alpha_i}(x)}{1-\bar{\alpha}_i\bar{F}(x)}\\&=&\frac{r(x)}{F(x)}\frac{\varphi'\left(\sum_{i=1}^n
\phi\left(\bar{F}_{\alpha_i}(x)\right)\right)}{\varphi\left(\sum_{i=1}^n
\phi\left(\bar{F}_{\alpha_i}(x)\right)\right)}\sum_{i=1}^n
\frac{\bar{F}_{\alpha_i}(x)}{\varphi'\left(\phi\left(\bar{F}_{\alpha_i}(x)\right)\right)}\frac{F(x)}{1-\bar{\alpha}_i\bar{F}(x)}\\
&=&\frac{r(x)}{F(x)}I_1\left(\phi\left(\bar{F}_{\alpha_1}(x)\right),\ldots,\phi\left(\bar{F}_{\alpha_n}(x)\right)\right),\end{eqnarray*}
where
\begin{eqnarray*}I_1\left(\phi\left(\bar{F}_{\alpha_1}(x)\right),\ldots,\phi\left(\bar{F}_{\alpha_n}(x)\right)\right)=\frac{\varphi'\left(\sum_{i=1}^n
\phi\left(\bar{F}_{\alpha_i}(x)\right)\right)}{\varphi\left(\sum_{i=1}^n
\phi\left(\bar{F}_{\alpha_i}(x)\right)\right)}
\sum_{i=1}^n\frac{\varphi\left(\phi\left(\bar{F}_{\alpha_i}(x)\right)\right)}{\varphi'\left(\phi\left(\bar{F}_{\alpha_i}(x)\right)\right)}\left(1-\varphi\left(\phi\left(\bar{F}_{\alpha_i}(x)\right)\right)\right).\end{eqnarray*}
It is easy to check that $\phi\left(\bar{F}_{\alpha_i}(x)\right)$ is
decreasing and convex in $\alpha_i$. From Theorem A.2 (Chapter 5) of
Marshall et al. \cite{marsh3}, $\boldsymbol\alpha
\stackrel{w}\succeq\boldsymbol\beta~\text{implies}~\left(\phi\left(\bar{F}_{\alpha_1}(x)\right),\ldots,\phi\left(\bar{F}_{\alpha_n}(x)\right)\right)\succeq_{w}
\left(\phi\left(\bar{F}_{\beta_1}(x)\right),\ldots,\phi\left(\bar{F}_{\beta_n}(x)\right)\right)$.
From Lemma \ref{lelu}, $I_1(\textbf{u})$ is increasing in $u_i$ for
$i\in I_n$ and Schur-convex with respect to $\textbf{u}$ whenever
$\varphi$ is log-concave and $\frac{\varphi (1-\varphi)}{\varphi'}$
is decreasing and concave. Then from Theorem A.8 (Chapter 3) of
Marshall et al. \cite{marsh3}, we get
$$I_1\left(\phi\left(\bar{F}_{\alpha_1}(x)\right),\ldots,\phi\left(\bar{F}_{\alpha_n}(x)\right)\right)\geq
I_1\left(\phi\left(\bar{F}_{\beta_1}(x)\right),\ldots,\phi\left(\bar{F}_{\beta_n}(x)\right)\right)$$
which implies $r_{X_{1:n}}(x)\geq r_{Y_{1:n}}(x)$, that is $X_{1:n}\leq_{hr}Y_{1:n}$.\\
Next we prove the theorem when $\frac{\varphi
(1-\varphi)}{\varphi'}$ is convex. Let
$z_i=\phi\left(\bar{F}_{\alpha_i}(x)\right)$. Then the hazard rate
function is given by
$$r_{X_{1:n}}(x)=\frac{r(x)}{F(x)}\frac{\varphi'\left(\sum_{i=1}^n
z_i\right)}{\varphi\left(\sum_{i=1}^n
z_i\right)}\sum_{i=1}^n\frac{\varphi\left(z_i\right)}{\varphi'\left(z_i\right)}\left(1-\varphi\left(z_i\right)\right).$$
Now, for $s\in I_n$,
\begin{eqnarray*}\frac{r_{X_{1:n}}(x)}{\partial \alpha_s}&=& \frac{r(x)}{F(x)} \left[\frac{\partial}{\partial z_s} \left(\frac{\varphi'\left(\sum_{i=1}^n z_i\right)}{\varphi\left(\sum_{i=1}^n
z_i\right)}\right)\frac{\partial z_s}{\partial \alpha_s}
\sum_{i=1}^n\frac{\varphi\left(z_i\right)\left(1-\varphi\left(z_i\right)\right)}{\varphi'\left(z_i\right)}\right.+\\&&\left.
\frac{\varphi'\left(\sum_{i=1}^n
z_i\right)}{\varphi\left(\sum_{i=1}^n z_i\right)}
\frac{\partial}{\partial z_s}
\left(\frac{\varphi\left(z_s\right)\left(1-\varphi\left(z_s\right)\right)}{\varphi'\left(z_s\right)}\right)
\frac{\partial z_s}{\partial \alpha_s}\right].
\end{eqnarray*}
Note that $z_s$ is decreasing in $\alpha_s$ and $\frac{\partial
z_s}{\partial \alpha_s}$ is increasing in $\alpha_s$. Since
$\varphi$ is log-concave and $\frac{\varphi (1-\varphi)}{\varphi'}$
is decreasing, we have $\frac{r_{X_{1:n}}(x)}{\partial \alpha_s}\leq
0$. Again
$$\frac{\partial}{\partial z_s} \left(\frac{\varphi'\left(\sum_{i=1}^n
z_i\right)}{\varphi\left(\sum_{i=1}^n
z_i\right)}\right)=\frac{\partial}{\partial z_t}
\left(\frac{\varphi'\left(\sum_{i=1}^n
z_i\right)}{\varphi\left(\sum_{i=1}^n
z_i\right)}\right),~\text{for}~s\neq t.$$ For $s\neq t$,
\begin{eqnarray*}&&(\alpha_s-\alpha_t)\left(\frac{r_{X_{1:n}}}{\partial
\alpha_s}-\frac{r_{X_{1:n}}}{\partial
\alpha_t}\right)\\&\stackrel{sign}{=}&(\alpha_s-\alpha_t)\left(\frac{\partial
z_s}{\partial \alpha_s}-\frac{\partial z_t}{\partial
\alpha_t}\right)+(\alpha_s-\alpha_t)\frac{\varphi'\left(\sum_{i=1}^n
z_i\right)}{\varphi\left(\sum_{i=1}^n z_i\right)}\times\\&&
\left[\left(-\frac{\partial}{\partial z_s}
\left(\frac{\varphi\left(z_s\right)\left(1-\varphi\left(z_s\right)\right)}{\varphi'\left(z_s\right)}\right)\right)\left(-\frac{\partial
z_s}{\partial \alpha_s}\right)-\left(-\frac{\partial}{\partial z_t}
\left(\frac{\varphi\left(z_t\right)\left(1-\varphi\left(z_t\right)\right)}{\varphi'\left(z_t\right)}\right)\right)\left(-\frac{\partial
z_t}{\partial \alpha_t}\right)\right]\\&\leq& 0,
\end{eqnarray*}
whenever $\frac{\varphi (1-\varphi)}{\varphi'}$ is convex in
addition to the log-concave $\varphi$ and decreasing $\frac{\varphi
(1-\varphi)}{\varphi'}$. Thus we have $r_{X_{1:n}}(x)$ is decreasing
in $\alpha_i$, $i\in I_n$ and Schur-convex in
$\boldsymbol\alpha=(\alpha_1,\alpha_2,...,\alpha_n)$. Then from
Lemma \ref{lewm}, we have $$\boldsymbol\alpha
\stackrel{w}{\succeq}\boldsymbol\beta~\text{implies}~r_{X_{1:n}}(x)\geq
r_{X_{1:n}}(x).$$ Hence the theorem follows.
\begin{coro}\label{corhr}Suppose the lifetime vectors $X\sim
PO(\bar{F},\boldsymbol\alpha,\varphi)$ and $Y\sim PO(\bar{F},\alpha
\boldsymbol 1,\varphi)$. Then, $X_{1:n}\leq_{hr}Y_{1:n}$ if $\alpha
\geq \frac{1}{n}\sum_{i=1}^n \alpha_i$, $\varphi$ is log-concave and
$\frac{\varphi (1-\varphi)}{\varphi'}$ is decreasing and concave (or
convex). This follows from the Theorem \ref{thhr} and using the fact
that
$(\alpha_1,\alpha_2,\ldots,\alpha_n)\stackrel{w}\succeq(\underbrace{\alpha,\alpha,\ldots,\alpha}_{n\;terms})$,
for $\alpha\geqslant\frac{1}{n}\sum_{i=1}^n\alpha_i$. $\hfill\Box$
\end{coro}
\begin{remark}\normalfont
It is to be noted that Archimedean copulas with generators
$\varphi(t)=2/(1+e^t)$ and $\varphi(t)=(-1+\theta)/(-e^t+\theta)$
for $-1\leq\theta\leq 0$ are some examples of survival copula such
that $\varphi$ is log-concave, and $\frac{\varphi
(1-\varphi)}{\varphi'}$ is decreasing and convex.
\end{remark}
The following counterexample shows that one may not get the the
ordering result in Theorem \ref{thhr} if the sufficient conditions
on the generator functions are dropped.
\begin{counterexample}\label{exthhr}\normalfont
        Consider two series systems, each
        comprising of three dependent and heterogeneous components with respective
        hazard rate functions $r_{X_{1:3}}(x)$ and
        $r_{Y_{1:3}}(x)$, with common baseline survival function
        $\bar{F}(x)=e^{-(0.5x)^2}$, $x\geq
        0$, $\boldsymbol\alpha=(0.2,0.4,0.6)$,
        $\boldsymbol\beta=(0.35,0.55,0.95)$ so that $\boldsymbol\alpha
        \stackrel{w}{\succeq}\boldsymbol\beta$. First we take the common generator $\varphi(x)=\log(e+x)^{-1/a}$, $a=0.1$, which is not log-concave
        but $\frac{\varphi (1-\varphi)}{\varphi'}$ is decreasing and convex. Next we take
        $\varphi(x)=(2/(1+e^x))^{1/a}$, $a=0.2$, which is log-concave but $\frac{\varphi (1-\varphi)}{\varphi'}$ is neither decreasing nor convex.
        For both the cases $r_{X_{1:3}}(x)$ and $r_{Y_{1:3}}(x)$ are depicted in Figure \ref{picth3.3p}(a) and \ref{picth3.3p}(b) respectively for some finite range of $x$.
        From both the figures we observe that the hazard rate ordering result in Theorem~\ref{thhr} is not attained.$\hfill\Box$
        \begin{figure}\begin{center}
                \includegraphics[width=15cm]{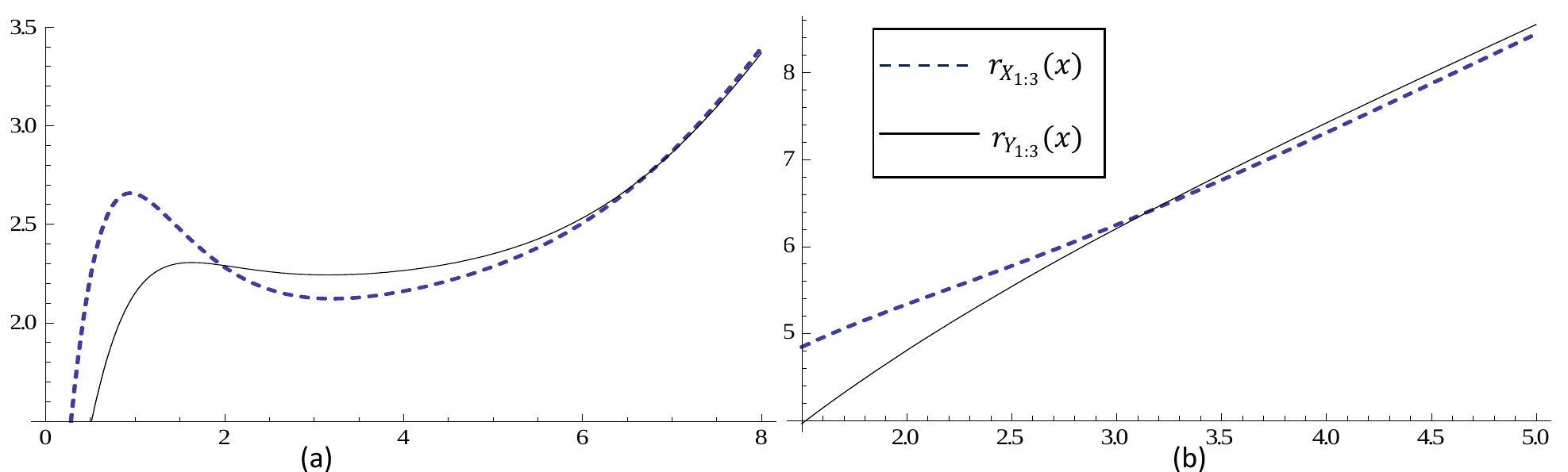}\\
                \caption{Plots of $r_{X_{1:3}}(x)$ and $r_{Y_{1:3}}(x)$ for (a) $\varphi(x)$ is not log-concave, (b) $\frac{\varphi (1-\varphi)}{\varphi'}$ is neither decreasing nor convex.}\label{picth3.3p}\end{center}
        \end{figure}
    \end{counterexample}
\section{Parallel systems with dependent and heterogeneous component lifetimes following PO Model}
Here, we compare the lifetimes of two parallel systems consisting of
dependent and heterogeneous components having lifetimes following
the PO model, with respect to some stochastic orders.\par Let the
lifetime vector $X=(X_1,X_2,...,X_n)$ be a set of dependent random
variables following the PO model with baseline survival function
$\bar{F}$ and having the joint distribution function coupled with
Archimedean survival copula with generator $\varphi$, denoted as
$X\sim PO(\bar{F},\boldsymbol\alpha,\varphi)$, where
$\boldsymbol\alpha=(\alpha_1,\alpha_2,...,\alpha_n)\in
\mathbb{R}^n_{+}$ is the proportional odds ratio vector. The
distribution function of $X_i$ is $F_{\alpha_i}(x)=\frac{
F(x)}{1-\bar{\alpha}_i\bar{F}(x)}$. The distribution function of
$X_{n:n}$ is given by
$$F_{X_{n:n}}(x)=P(X_{n:n}\leq x)=P(X_{i}<x,i\in
I_n)=\varphi\left(\sum_{i=1}^n
\phi\left(F_{\alpha_i}(x)\right)\right)=S_2(F(x),\boldsymbol\alpha,\varphi),~say,$$
where $\phi(u)=\varphi^{-1}(u)$, $u\in(0,1]$. The reversed hazard
rate function of $X_{n:n}$ is obtained as
\begin{equation}\label{prhr}
\tilde{r}_{X_{n:n}}(x)=\tilde{r}(x) \frac{\varphi'\left(\sum_{i=1}^n
\phi\left(F_{\alpha_i}(x)\right)\right)}{\varphi\left(\sum_{i=1}^n
\phi\left(F_{\alpha_i}(x)\right)\right)}\sum_{i=1}^n
\phi'\left(F_{\alpha_i}(x)\right)\frac{\alpha_i
F_{\alpha_i}(x)}{1-\bar{\alpha}_i \bar{F}(x)},
\end{equation} where $\tilde{r}$ denotes the baseline reversed hazard rate
function.
\begin{lem}\label{lej2}
For any $x\in[0,1]$, $S_2(x,\boldsymbol\alpha,\varphi)$ is
decreasing in $\alpha_i$, $i\in I_n$. Furthermore $S_2$ is
Schur-convex with respect to $\boldsymbol\alpha$ whenever $\varphi$
is log-concave.
 \end {lem}
\textbf{Proof:} For $s\in I_n$,
$$\frac{\partial S_2}{\partial \alpha_s}=-\varphi'\left(\sum_{i=1}^n
\phi\left(\frac{x}{1-\bar{\alpha}_i
(1-x)}\right)\right)\phi'\left(\frac{x}{1-\bar{\alpha}_s
(1-x)}\right)\frac{x(1-x)}{(1-\bar{\alpha}_s (1-x))^2}.$$ Since both
$\varphi(u)$ and $\phi(u)$ are decreasing for all $u\geq 0$,
$\frac{\partial S_2}{\partial \alpha_s}\leq 0$. So
$S_2(x,\boldsymbol\alpha,\varphi)$ is decreasing in $\alpha_i$ for
any $x\in[0,1]$.\\For $s\neq t$,
\begin{eqnarray*}&&(\alpha_s-\alpha_t)\left(\frac{\partial S_2}{\partial
\alpha_s}-\frac{\partial S_2}{\partial
\alpha_t}\right)\\&=&-(\alpha_s-\alpha_t)\varphi'\left(\sum_{i=1}^n
\phi\left(v_i\right)\right)\left[\phi'\left(v_s\right)\frac{x(1-x)}{(1-\bar{\alpha}_s
x)^2}-\phi'\left(v_t\right)\frac{x(1-x)}{(1-\bar{\alpha}_t x)^2}\right],~v_i=\frac{x}{1-\bar{\alpha}_i (1-x)}\\
&\stackrel{sign}{=}&(\alpha_s-\alpha_t)\left[-\left(-\frac{\varphi(\phi(v_s))}{\varphi'(\phi(v_s))}\right)\frac{1}{1-\bar{\alpha}_s
x}+\left(-\frac{\varphi(\phi(v_t))}{\varphi'(\phi(v_t))}\right)\frac{1}{1-\bar{\alpha}_t
x}\right]
\\&\geq& 0,\end{eqnarray*}
where the last inequality is derived using the fact that $\varphi$
is log-concave. So from Lemma \ref{lesc}, $S_2$ is Schur-convex in
$\boldsymbol\alpha=(\alpha_1,\alpha_2,...,\alpha_n)$.$\hfill\Box$

Suppose there are two parallel systems with lifetime vectors
$X=(X_1,X_2,...,X_n)$ and $Y=(Y_1,Y_2,...,Y_n)$, formed out of $n$
dependent and heterogeneous components where the component lifetimes
follow PO model. The following theorem compares the lifetimes of two
such parallel systems in the sense of usual stochastic order.
\begin{thm}\label{thstp}Suppose the lifetime vectors $X\sim PO(\bar{F},\boldsymbol\alpha,\varphi_1)$ and $Y\sim
PO(\bar{F},\boldsymbol\beta,\varphi_2)$. If $\varphi_1$ or
$\varphi_2$ is log-concave and $\phi_1\circ\varphi_2$ is
superadditive, then
$$\boldsymbol\alpha \stackrel{w}{\succeq}\boldsymbol\beta~\text{implies}~ X_{n:n}\leq_{st}Y_{n:n}.$$
\end{thm}
\textbf{Proof:} If $\varphi_1$ is log-concave, then from Lemma
\ref{lej2} and Lemma \ref{lewm}, we have
\begin{equation}\label{pst}\boldsymbol\alpha
\stackrel{w}{\succeq}\boldsymbol\beta~\text{implies}~S_2(F(x),\boldsymbol\alpha,\varphi_1)\geq
S_2(F(x),\boldsymbol\beta,\varphi_1).\end{equation} Since
$\phi_1\circ\varphi_2$ is superadditive, so from Lemma \ref{lecop}
(by replacing $\varphi_1$ by $\varphi_2$ and vice versa), we have
\begin{equation}\label{pst2} S_2(F(x),\boldsymbol\beta,\varphi_1)\geq
S_2(F(x),\boldsymbol\beta,\varphi_2).\end{equation} Combining
(\ref{pst}) and (\ref{pst2}), we get
$S_2(F(x),\boldsymbol\alpha,\varphi_1)\geq
S_2(F(x),\boldsymbol\beta,\varphi_2)$. That is
$X_{n:n}\leq_{st}Y_{n:n}$. \\ Now suppose $\varphi_2$ is
log-concave, then
\begin{eqnarray*}S_2(F(x),\boldsymbol\alpha,\varphi_1)&\geq&
S_2(F(x),\boldsymbol\alpha,\varphi_2)\\&\geq&
S_2(F(x),\boldsymbol\beta,\varphi_2),\end{eqnarray*} where the first
inequality follows from the fact that $\phi_1\circ\varphi_2$ is
superadditive, whereas the second inequality follows from the fact
that $\boldsymbol\alpha \stackrel{w}{\succeq}\boldsymbol\beta$. This
proves the result.$\hfill\Box$
\begin{coro}\label{corstp}Suppose the lifetime vectors $X\sim PO(\bar{F},\boldsymbol\alpha,\varphi)$ and $Y\sim
PO(\bar{F},\boldsymbol\beta,\varphi)$. If $\varphi$ is log-concave,
then
$$\boldsymbol\alpha \stackrel{w}{\succeq}\boldsymbol\beta~\text{implies}~
X_{n:n}\leq_{st}Y_{n:n}.\hfill\Box$$
\end{coro}
The following counterexample shows that one may not get the the
ordering result in Theorem \ref{thstp} if the sufficient conditions
on the generator functions are dropped.
\begin{counterexample}\label{exthstp}\normalfont
        Consider two parallel systems, each
        comprising of three dependent and heterogeneous components with respective
        distribution functions $F_{X_{3:3}}(x)=S_2(F(x),\boldsymbol\alpha,\varphi_1)$ and
        $F_{Y_{3:3}}(x)=S_2(F(x),\boldsymbol\beta,\varphi_2)$, where
        $F(x)=1-e^{-x^{0.5}}$, $x\geq
        0$, $\boldsymbol\alpha=(0.9,1.45,2.15)$,
        $\boldsymbol\beta=(1.2,1.95,2.65)$ so that $\boldsymbol\alpha
        \stackrel{w}{\succeq}\boldsymbol\beta$. First we take $\varphi_1(x)=\theta_1/\log(x+e^{\theta_1})$ and
        $\varphi_2(x)=e^{1-(1+x)^{1/\theta_2}}$ with $\theta_1=0.9$ and $\theta_2=8$ so that neither $\varphi_1$ nor $\varphi_2$ is log-concave
        but $\phi_1\circ \varphi_2$ is super additive. Next we take $\varphi_1(x)=e^{(1-e^x)/\theta_1}$ and
        $\varphi_2(x)=(2/(e^x+1))^{1/\theta_2}$ with $\theta_1=0.9$ and $\theta_2=0.2$ so that $\varphi_1$ is log-concave
        but $\phi_1\circ \varphi_2$ is not super additive. For both the cases $F_{X_{3:3}}(x)$ and $F_{Y_{3:3}}(x)$ are depicted in Figure \ref{picth4.1p}(a) and \ref{picth4.1p}(b) respectively for some finite range of $x$.
        From both the figures we observe that the stochastic ordering result in Theorem~\ref{thstp} is not attained.$\hfill\Box$
        \begin{figure}\begin{center}
                \includegraphics[width=14cm]{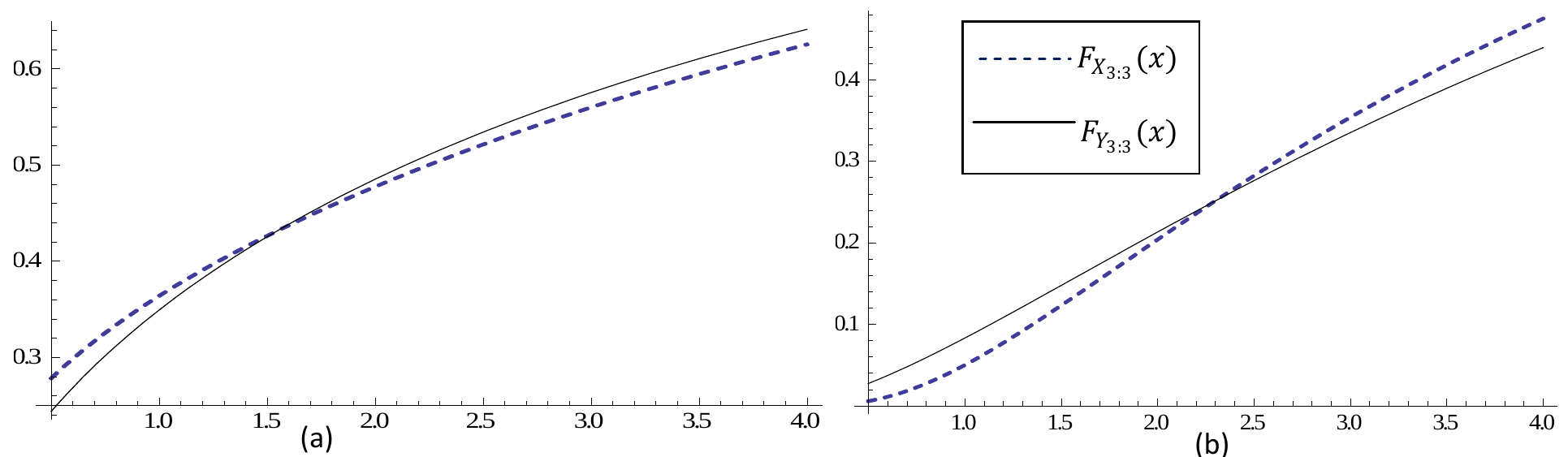}\\
                \caption{Plots of $F_{X_{3:3}}(x)$ and $F_{Y_{3:3}}(x)$ for (a) neither $\varphi_1$ nor $\varphi_2$ is log-concave, (b) $\phi_1\circ \varphi_2$ is not super additive.}\label{picth4.1p}\end{center}
        \end{figure}
    \end{counterexample}
Next we show the reversed hazard rate order of lifetimes of two
parallel systems of dependent and heterogeneous components.
\begin{thm}\label{thrhr}Suppose the lifetime vectors $X\sim PO(\bar{F},\boldsymbol\alpha,\varphi)$ and $Y\sim
PO(\bar{F},\boldsymbol\beta,\varphi)$. If $\varphi$ is log-concave
and $\frac{\varphi (1-\varphi)}{\varphi'}$ is decreasing and convex,
then
$$\boldsymbol\alpha \stackrel{w}\succeq\boldsymbol\beta~\text{implies}~ X_{n:n}\leq_{rhr}Y_{n:n}.$$
\end{thm}
\textbf{Proof:} From (\ref{prhr}), the reversed hazard function of
$X_{n:n}$ is given by
\begin{eqnarray*}\tilde{r}_{X_{n:n}}(x)&=&\frac{\tilde{r}(x)}{\bar{F}(x)} \frac{\varphi'\left(\sum_{i=1}^n
\phi\left(F_{\alpha_i}(x)\right)\right)}{\varphi\left(\sum_{i=1}^n
\phi\left(F_{\alpha_i}(x)\right)\right)}\sum_{i=1}^n
\frac{F_{\alpha_i}(x)}{\varphi'\left(\phi\left(F_{\alpha_i}(x)\right)\right)}\bar{F}_{\alpha_i}(x)
\\&=&\frac{\tilde{r}(x)}{\bar{F}(x)}
\frac{\varphi'\left(\sum_{i=1}^n
\phi\left(F_{\alpha_i}(x)\right)\right)}{\varphi\left(\sum_{i=1}^n
\phi\left(F_{\alpha_i}(x)\right)\right)}
\sum_{i=1}^n\frac{\varphi\left(\phi\left(F_{\alpha_i}(x)\right)\right)}{\varphi'\left(\phi\left(F_{\alpha_i}(x)\right)\right)}\left(1-\varphi\left(\phi\left(F_{\alpha_i}(x)\right)\right)\right)\\
&=&\frac{\tilde{r}(x)}{\bar{F}(x)} \frac{\varphi'\left(\sum_{i=1}^n
\xi_i\right)}{\varphi\left(\sum_{i=1}^n
\xi_i\right)}\sum_{i=1}^n\frac{\varphi\left(\xi_i\right)}{\varphi'\left(\xi_i\right)}\left(1-\varphi\left(\xi_i\right)\right),\end{eqnarray*}
where $\xi_i=\phi\left(F_{\alpha_i}(x)\right)$. Now, for $s\in I_n$,
\begin{eqnarray*}\frac{\tilde{r}_{X_{n:n}}(x)}{\partial \alpha_s}&=& \frac{\tilde{r}(x)}{\bar{F}(x)} \left[\frac{\partial}{\partial \xi_s} \left(\frac{\varphi'\left(\sum_{i=1}^n \xi_i\right)}{\varphi\left(\sum_{i=1}^n
\xi_i\right)}\right)\frac{\partial \xi_s}{\partial \alpha_s}
\sum_{i=1}^n\frac{\varphi\left(\xi_i\right)\left(1-\varphi\left(\xi_i\right)\right)}{\varphi'\left(\xi_i\right)}\right.+\\&&\left.
\frac{\varphi'\left(\sum_{i=1}^n
\xi_i\right)}{\varphi\left(\sum_{i=1}^n \xi_i\right)}
\frac{\partial}{\partial \xi_s}
\left(\frac{\varphi\left(\xi_s\right)\left(1-\varphi\left(\xi_s\right)\right)}{\varphi'\left(\xi_s\right)}\right)
\frac{\partial \xi_s}{\partial \alpha_s}\right].
\end{eqnarray*}
Note that $\xi_s$ is increasing in $\alpha_s$ and $\frac{\partial
\xi_s}{\partial \alpha_s}$ is decreasing in $\alpha_s$. Since
$\varphi$ is log-concave and $\frac{\varphi (1-\varphi)}{\varphi'}$
is decreasing, we have $\frac{\tilde{r}_{X_{n:n}}(x)}{\partial
\alpha_s}\geq 0$. Again
$$\frac{\partial}{\partial \xi_s} \left(\frac{\varphi'\left(\sum_{i=1}^n
\xi_i\right)}{\varphi\left(\sum_{i=1}^n
\xi_i\right)}\right)=\frac{\partial}{\partial \xi_t}
\left(\frac{\varphi'\left(\sum_{i=1}^n
\xi_i\right)}{\varphi\left(\sum_{i=1}^n
\xi_i\right)}\right),~\text{for}~s\neq t.$$ For $s\neq t$,
\begin{eqnarray*}&&(\alpha_s-\alpha_t)\left(\frac{\tilde{r}_{X_{n:n}}}{\partial
\alpha_s}-\frac{\tilde{r}_{X_{n:n}}}{\partial
\alpha_t}\right)\\&\stackrel{sign}{=}&(\alpha_s-\alpha_t)\left(\frac{\partial
\xi_s}{\partial \alpha_s}-\frac{\partial \xi_t}{\partial
\alpha_t}\right)+(\alpha_s-\alpha_t)\left(-\frac{\varphi'\left(\sum_{i=1}^n
\xi_i\right)}{\varphi\left(\sum_{i=1}^n
\xi_i\right)}\right)\times\\&& \left[\left(-\frac{\partial}{\partial
\xi_s}
\left(\frac{\varphi\left(\xi_s\right)\left(1-\varphi\left(\xi_s\right)\right)}{\varphi'\left(\xi_s\right)}\right)\right)\frac{\partial
\xi_s}{\partial \alpha_s}-\left(-\frac{\partial}{\partial \xi_t}
\left(\frac{\varphi\left(\xi_t\right)\left(1-\varphi\left(\xi_t\right)\right)}{\varphi'\left(\xi_t\right)}\right)\right)\frac{\partial
\xi_t}{\partial \alpha_t}\right]\\&\leq& 0,
\end{eqnarray*}
as $\frac{\varphi (1-\varphi)}{\varphi'}$ is decreasing and convex.
Thus we have $\tilde{r}_{X_{n:n}}(x)$ is increasing in $\alpha_i$,
$i\in I_n$ and Schur-concave in
$\boldsymbol\alpha=(\alpha_1,\alpha_2,...,\alpha_n)$. Then from
Lemma \ref{lewm}, we have $$\boldsymbol\alpha
\stackrel{w}{\succeq}\boldsymbol\beta~\text{implies}~\tilde{r}_{X_{n:n}}(x)\leq
\tilde{r}_{Y_{n:n}}(x).$$ Hence the theorem follows.
The following counterexample shows that one may not get the the
ordering result in Theorem \ref{thrhr} if the sufficient conditions
on the generator functions are dropped.
\begin{counterexample}\label{exthrhr}\normalfont
        Consider two parallel systems, each
        comprising of four dependent and heterogeneous components with respective
        reversed hazard rate functions $\tilde{r}_{X_{4:4}}(x)$ and
        $\tilde{r}_{Y_{4:4}}(x)$, with common baseline survival function
        $\bar{F}(x)=e^{-x^3}$, $x\geq
        0$, $\boldsymbol\alpha=(0.2,0.6,1.5,3.5)$,
        $\boldsymbol\beta=(0.8,0.9,4.5,5.5)$ so that $\boldsymbol\alpha
        \stackrel{w}{\succeq}\boldsymbol\beta$.
        First we take the common generator $\varphi(x)=(1/(ax+1))^{1/a}$, $a=0.2$, which is not log-concave
        but $\frac{\varphi (1-\varphi)}{\varphi'}$ is decreasing and convex. Next we take
        $\varphi(x)=(2/(1+e^x))^{1/a}$, $a=0.2$, which is log-concave but $\frac{\varphi (1-\varphi)}{\varphi'}$ is neither decreasing nor convex.
        For both the cases $\tilde{r}_{X_{4:4}}(x)$ and $\tilde{r}_{Y_{4:4}}(x)$ are depicted in Figure \ref{picth4.2p}(a) and \ref{picth4.2p}(b) respectively for some finite range of $x$.
        From both the figures we observe that the reversed hazard rate ordering result in Theorem~\ref{thrhr} is not attained.$\hfill\Box$
        \begin{figure}\begin{center}
                \includegraphics[width=15cm]{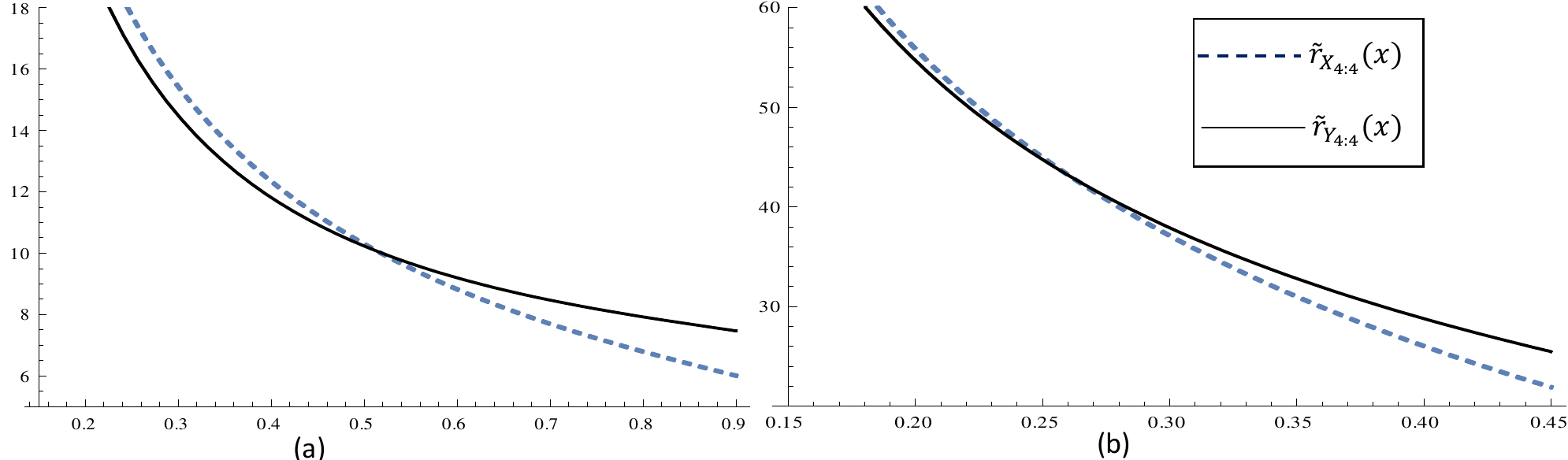}\\
                \caption{Plots of $\tilde{r}_{X_{4:4}}(x)$ and $\tilde{r}_{Y_{4:4}}(x)$ for (a) $\varphi(x)$ is not log-concave, (b) $\frac{\varphi (1-\varphi)}{\varphi'}$ is neither decreasing nor convex.}\label{picth4.2p}\end{center}
        \end{figure}
    \end{counterexample}
\section{Applications}
In this section, we highlight some potential applications of the
established results. Consider a series (or a parallel) system of $n$
components having dependent lifetimes. It is quite practical that
the odds functions of all the components (i.e. the odds of surviving
beyond a specified time $t$) may not be the same for various
possible reasons, like the components are manufactured by different
manufacturing units or they are subjected to different levels of
stress. So let the odds function of $i$th component is proportional
to a baseline odds function with proportionality constant (odds
ratio) $\alpha_i$, $i=1,2,\ldots,n$. Now consider another series (or
parallel) system of $n$ dependent components having different odds
ratios $\beta_i$, $i=1,2,\ldots,n$. Even if a same system operates
in two different levels of environments/stress (e.g., voltage,
temperature, compression and tension), then reliability
characteristics (e.g., odds function) of a component of the system
generally will not be the same in the two different environments. So
it is a subject of interest to compare lifetimes of two such
systems, i.e. under what conditions one system will be more reliable
than other. Theorems \ref{thst} and \ref{thst1} (resp. Theorem
\ref{thstp}) give the conditions on the corresponding odds ratio
vectors and the generators of the survival copulas under which a
series (resp. parallel) system will have stochastically longer
lifetime than that of the other. Similarly Theorem \ref{thhr} (resp.
Theorem \ref{thrhr}) gives the conditions under which failure
rate of a series (resp. parallel) system will be smaller than that
of the other.

Next we show that using our proposed results one can compare the
lifetime of two series systems whose components are subjected to
random shock instantaneously \cite{fangL}. Suppose random variable
$X_i$ denotes the lifetime of $i$-th component of the series system.
Define Bernoulli random variable $I_{p_i}$ associated with $X_i$,
where $I_{p_i}=1$ if shock does not occur and 0 if shock occurs with
$p_i=P\left(I_{p_i}=1\right)$, $i=1,\ldots,n$. Assume that
$I_{p_1},\ldots,I_{p_n}$ are independent random variables, and also
they are independent of $X_1,\ldots,X_n$. Let $X_i^{*}=X_i I_{p_i}$,
$i=1,\ldots,n$, and denote $X_{1:n}^*=\min(X_1^{*},\ldots,X_n^{*})$.
Similarly assume that $I_{q_1},\ldots,I_{q_n}$ are independent
Bernoulli random variables, and also they are independent of $Y_i$'s
with $q_i=P\left(I_{q_i}=1\right)$, $i=1,\ldots,n$. Denote
$Y_{1:n}^*=\min(Y_1^{*},\ldots,Y_n^{*})$, where $Y_i^{*}=Y_i
I_{q_i}$, $i=1,\ldots,n$. Here $X_{1:n}^*$ represents the lifetime
of a series system whose components are subjected to random shock
instantaneously. Similarly $Y_{1:n}^*$ represents the the lifetime
of another such series system. Now, if $X\sim
PO(\bar{F},\boldsymbol\alpha,\varphi_1)$ and $Y\sim
PO(\bar{F},\boldsymbol\beta,\varphi_2)$, then with the help of the
Theorems \ref{thst}, \ref{thst1}, \ref{thhr}, and the associated
corollaries \ref{corst}, \ref{corst1}, we can establish following
stochastic comparisons between such smallest order statistics from
the fact that $P(X_{1:n}^*>x)=P(X_1>x,\ldots,X_n>x)P(I_{p_i}=1,i\in
I_n)=P(X_{1:n}>x)\prod_{i}^n p_i$.
\begin{thm}\label{thstca}Suppose the lifetime vectors $X\sim PO(\bar{F},\boldsymbol\alpha,\varphi_1)$ and $Y\sim
PO(\bar{F},\boldsymbol\beta,\varphi_2)$. If $\varphi_1$ or
$\varphi_2$ is log-convex, $\phi_2\circ\varphi_1$ is superadditive
and $\prod_{i}^n p_i\leq \prod_{i}^n q_i$, then
$$\boldsymbol\alpha \stackrel{p}{\succeq}\boldsymbol\beta~\text{implies}~ X_{1:n}^*\leq_{st}Y_{1:n}^*.$$
\end{thm}
\begin{coro}\label{corstca}Suppose the lifetime vectors $X\sim PO(\bar{F},\boldsymbol\alpha,\varphi)$ and $Y\sim
PO(\bar{F},\boldsymbol\beta,\varphi)$. If $\varphi$ is log-convex
and $\prod_{i}^n p_i\leq \prod_{i}^n q_i$, Then
$$\boldsymbol\alpha \stackrel{p}{\succeq}\boldsymbol\beta~\text{implies}~
X_{1:n}^*\leq_{st}Y_{1:n}^*.$$
\end{coro}
\begin{thm}\label{thst1ca}Suppose the lifetime vectors $X\sim PO(\bar{F},\boldsymbol\alpha,\varphi_1)$ and $Y\sim
PO(\bar{F},\boldsymbol\beta,\varphi_2)$. If $\phi_2\circ\varphi_1$
is superadditive and $\prod_{i}^n p_i\leq \prod_{i}^n q_i$, then
$$\boldsymbol\alpha \stackrel{w}{\succeq}\boldsymbol\beta~\text{implies}~ X_{1:n}^*\leq_{st}Y_{1:n}^*.$$
\end{thm}
\begin{coro}\label{corst1ca}Suppose the lifetime vectors $X\sim PO(\bar{F},\boldsymbol\alpha,\varphi)$ and $Y\sim
PO(\bar{F},\boldsymbol\beta,\varphi)$. If $\prod_{i}^n p_i\leq
\prod_{i}^n q_i$, then
$$\boldsymbol\alpha \stackrel{w}{\succeq}\boldsymbol\beta~\text{implies}~
X_{1:n}^*\leq_{st}Y_{1:n}^*.$$
\end{coro}
\begin{thm}\label{thhrca}Suppose the lifetime vectors $X\sim PO(\bar{F},\boldsymbol\alpha,\varphi)$ and $Y\sim
PO(\bar{F},\boldsymbol\beta,\varphi)$. If $\varphi$ is log-concave
and $\frac{\varphi (1-\varphi)}{\varphi'}$ is decreasing and
concave, then
$$\boldsymbol\alpha \stackrel{w}\succeq\boldsymbol\beta~\text{implies}~ X_{1:n}^*\leq_{hr}Y_{1:n}^*.$$
\end{thm}
We will end this section by mentioning an other potential
application. In actuarial science, $X_{1:n}^*$ corresponds to the
smallest claim amount in a portfolio of risks \cite{barm,lil,zhang},
where $X_i$'s represent sizes of random claims of multiple risks
covered by a policy that can be made in an insurance period and the
corresponding $I_{p_i}$'s indicate the occurrence of these claims.
That means $I_{p_i}=1$ whenever the $i$th policy makes random claim
$X_i$ and $I_{p_i}=0$ whenever there is no claim. Similarly suppose
$Y_{1:n}^*$ represents the smallest claim amount in an another
portfolio of risks. The above theorems can be used in stochastic
comparisons between the smallest claim amounts of two different
portfolio of risks.

\section*{Acknowledgements:} The authors are thankful to the Area Editor and the anonymous Reviewers for the constructive comments and suggestions, which lead to an improved version of the manuscript.

\end{document}